\documentclass[12pt,leqno]{article}
\usepackage{amstext}
\usepackage{theorem}
\usepackage{amssymb}
\usepackage{latexsym}
\newfont{\ninemsbm}{msbm10 scaled 0900}
\newfont{\tenmsbm}{msbm10 scaled 1100}
\newfont{\nineeufb}{eufb10 scaled 0900}
\newfont{\teneufb}{eufb10 scaled 1100}
\newfont{\teneusm}{eusm10 scaled 1100}
\newfont{\nineeusm}{eusm10 scaled 0900}
\newfont{\paren}{cmex10 scaled 1100}

\newcommand{\sbm}[1]{\mbox{{\tenmsbm {#1}}}}

\newcommand{\fasm}[1]{\mbox{{\teneusm {#1}}}}

\newcommand{\lra}{\longrightarrow}

\newcommand{\be}{\begin{enumerate}}
\newcommand{\ee}{\end{enumerate}}

\newcommand{\br}{\begin{array}}
\newcommand{\er}{\end{array}}

\newcommand{\gchi}{\raisebox{.4ex}{$\chi$}}

\newcommand{\up}{^{\prime}}

\newcommand{\upp}{^{\prime\prime}}

\newcommand{\putref}[2]{
\noindent\begin{minipage}[t]{.07\linewidth}{#1}\end{minipage}
\begin{minipage}[t]{.9\linewidth}{#2}\end{minipage}}

\newcommand{\binm}[2]{\left(
\!\!\br{c}#1\vspace{-.5ex}\\#2\er\!\!\right)}
\newtheorem{thm}{Theorem.}[section]
\newtheorem{prop}[thm]{Proposition.}
\newtheorem{cor}[thm]{Corollary.}
\newtheorem{lem}[thm]{Lemma.}
\newtheorem{lemb}{Lemma.}[thm]

\newtheorem{rems}[thm]{Remarks.}
\newtheorem{defn}[thm]{Definition.}

\newtheorem{prop-def}[thm]{Proposition-Definition.}
\newtheorem{lemd}[thm]{Lemma-definition.}
\def\pmb#1{\setbox0=\hbox{#1}%
\kern-.025em\copy0\kern-\wd0
\kern.05em\copy0\kern-\wd0
\kern-.025em\raise.0433em\box0}

\begin{document}

\begin{center}{\Large\bf GENERIC HYPERSURFACE
SINGULARITIES}\end{center}

\noindent{Summary}

\noindent {\em The problem considered here can be viewed as the analogue in
higher dimensions of the one variable polynomial interpolation of Lagrange
and Newton. Let $x_1,\ldots ,x_r$ be closed points in general position in
projective space $\sbm P^n$, then the linear subspace $V$ of $H^0(\sbm
P^n,\fasm O(d))$ (the space of homogeneous polynomials of degree $d$ on
$\sbm P^n$ ) formed by those polynomials which are singular at each $x_i$,
is given by
$r(n+1)$ linear equations in the coefficients, expressing the fact that the
polynomial vanishes with its first derivatives at $x_1,\ldots ,x_r$. As such,
the ``expected'' value for the dimension of
$V$ is max$(0,h^0(\fasm O(d))-r(n+1))$. We prove that $V$ has the
``expected'' dimension for $d\geq 5$ (theorem A.). This theorem was first
proven in [A] using a very complicated induction with many initial cases.
Here we  give a greatly simplified proof using techniques developped by the
authors while treating the corresponding problem in lower degrees.}

\section{Introduction}

\noindent{\bf Theorem A.} {\em Let $x_1,\ldots,x_r $ be closed points in
general position in n dimensional projective space $\sbm P^n$ over an
algebraically closed field
$k$. Then the vector subspace of $H^0(\sbm P^n,\fasm O(d))$ of homogeneous
polynomials of degree $d\geq 5$ on $\sbm P^n$ having singularites 
at each $x_i$ has the ``expected'' dimension max$(0,h^0(\fasm
O(d))-r(n+1))$.}
	
We will give a rough outline of the proof. By a standard arguement we reduce
to a proposition (2.6) concerning unions of double points (i.e. first
infinitesimal neighbourhoods of closed points in $\sbm P^n$). This is done in
\S2. We attempt to prove this proposition by specialising a number of the
double points to have support in a hyperplane $H$. Each of the specialised
double points then induces a double point in $H$ (its trace on $H$) and has a
residual simple point. If the numbers fit (see `adjusting' 3.2) we obtain a
(lower dimensional) trace on $H$ of the kind treated by our proposition and
we can then argue by induction.

Unfortunately these numbers do not always fit! 
Nonetheless, we are able to overcome the problem by using the differential
lemma of [AH1] (see lemma 3.3). We begin as above by specialising a number
of double points to have support in $H$. We put in just enough double points
to have at least the required number of conditions in $H$. Any excess of
conditions in $H$, is then considered as coming from the trace of a single
double point $Z$ with support in $H$. The lemma then allows us to use the
induction hypothesis (applied in $H$) to reduce to a residual problem (see
degue in 3.3) in degree one less, involving a subscheme $d^{\prime}Z$ of
$Z\cap H$. This $d^{\prime}Z$ depends on the union of all the residual points
of the other double points we specialised to $H$ and although it is well
determined (see 3.4 (2)), it is very difficult to locate. In [AH1] we opted for
a strategy of specialisation, allowing us to identify $d^{\prime}Z$ at the
cost of a more complicated arguement. Here, we are able to proceed without
identifying $d^{\prime}Z$ (this means that the degue of 3.3 is true for all
subschemes of $Z\cap H$ having the same length as $d^{\prime}Z$). We do
this in 5.1, by ejecting (see 3.2) $d^{\prime}Z$ out of $H$ and specialising
further double points to have support in $H$. This delays the influence of
$d^{\prime}Z$ to a later stage in the argument, when we dispose of many
points in $H$ on whose position, $d^{\prime}Z$ does not depend. By
specialisation, we transform the union of $d^{\prime}Z$ with the right
number of simple points, into a double point of $\sbm P^n$. This leads to a
statement (6.1) in which there is no unkown component. The proof of 6.1,
uses simpler techniques of specialisation and ejection. 
	
	It is interesting to note that there are numerous examples in degrees
$d\leq5$ where  theorem A. fails. In fact for $d=2$ all values of $r$ with
$2\leq r\leq n$ give such examples (see [A] \S 4.), while for $d=4$, there is
an example for each $n= 2,3,4$ (see [A] \S 5.and [AH2]). The only example
with
$d=3$ is that of a hypercubic in $\sbm P^4$ with seven generic singularities
(see [A] \S 5. [CH] and [AH3]). Of course, unlike the one variable case, the
theorem fails miserably in higher dimensions if one abandons the general
position hypothesis on the $x_i$.
	
From another point of view, one can consider the general problem  of 
determining the dimension of the space of homogeneous polynomials of
degree d, having singularites of multiplicity $\geq m_i$ at $x_i$ for some
fixed sequence ($m_1,\ldots ,m_r$). For $n=2$, this has been persued further
in [H3], where the second author has formulated a conjecture giving a clear
geometric meaning to those cases where the ``expected'' dimension is not
realised. Eventually, theorem A. should be completed by a similar
interpretation of the examples in the preceeding paragraph.

\section{ A reformulation of the theorem}

Throughout the rest of the article we work over a 
fixed algebraically closed base field k. We write $\sbm P^n$ for 
projective $n$ space over this field.

In this section we formulate clearly theorem A in terms of projective
cohomology. Once this has been done we show that theorem A results from
proposition 1.6 below. The remaining sections are devoted to proving 1.6.
	 
	Let $x_1,\ldots , x_r$ be closed points in general position in  $\sbm P^n$
and let $Y$ be the union of the first infinitesimal neighbourhoods of the
$x_i$ ($i=1,\ldots ,r$). In any degree $d$, the homogeneous polynomials of
degree
$d$ which are singular at each $x_i$ are canonically identified, via the usual
exact sequence 
$$0 \lra H^0(\sbm P^n,\fasm I_Y(d))\lra  H^0(\sbm P^n,\fasm O(d))\lra
H^0(\sbm P^n,\fasm O_Y(d))\lra  H^1(\sbm P^n,\fasm I_Y(d))\lra 0$$
with $H^0(\sbm P^n,\fasm I_Y(d))$, where $\fasm I_Y$ is the ideal sheaf of
$Y$ as a closed subscheme of $\sbm P^n$. 

To say that $H^0(\sbm P^n, I_Y(d))$
has the expected dimension is equivalent to saying that the canonical map
$ H^0(\sbm P^n,\fasm O(d))\lra
H^0(\sbm P^n,\fasm O_Y(d))$ has maximal rank as a map of vector spaces or,
what amounts to the same thing, that one or the other of the two numbers
$h^0(\sbm P^n,\fasm I_Y(d))$ or $h^0(\sbm P^n,\fasm I_Y(d))$ is zero. This
motivates the

\begin{defn} Let $Y$ be a closed subscheme of $\sbm P^n$. We say that $Y$
has maximal rank in degree $d$ if the canonical map $H^0(\sbm P^n,\fasm
O(d))\lra H^0(\sbm P^n,\fasm O_Y(d))$
has maximal rank as a map of vector spaces. We say that $Y$ has maximal
rank if it has maximal rank in all degrees $d\geq 0$. 

We say that $Y$ is numerically adjusted (resp. adjusted) in degree $d$ if
$\gchi (\sbm P^n ,\fasm I_Y(d))=0$ (resp. $h^i(\sbm P^n , I_Y(d))=0$  for
$i\geq 0$).
\end{defn}

Clearly if $Y$ is adjusted in degree $d$ then it is both numerically adjusted
and of maximal rank in degree $d$.

The following proposition and its corollary show how these 
various notions are related for zero dimensional subschemes.

\begin{prop}.  Let $Y^{\prime}\subset Y\subset Y^{\prime\prime}$ be zero
dimensional closed subschemes of $\sbm P^n$. Fix $d\geq 0$ and consider the
following two maps
$$\br{llll}
		(i) 	&H^0(\sbm P^n,\fasm O(d^{\prime}))\lra  H^0(\sbm P^n,\fasm
O_{Y\up}(d^{\prime}))& ; &   d^{\prime}\geq 0\\
		(ii) &	H^0(\sbm P^n,\fasm O(d^{\prime\prime}))\lra H^0(\sbm P^n,\fasm
O_{Y^{\prime\prime}}(d^{\prime\prime}))& ;    &d^{\prime\prime}\geq 0
 \er$$
Then we have 
\newcounter{mist}
\begin{list}{\textbf{ (\alph{mist}) }}
{\setlength{\leftmargin}{1cm}
\setlength{\rightmargin}{1cm}
\usecounter{mist}  } 
\item If (i) is injective, then $H^0(\sbm P^n,\fasm O(d))\lra
 H^0(\sbm P^n,\fasm O_{Y}(d))$ is injective for $d\leq d^{\prime}$.
\item	If (ii) is surjective, then $H^0(\sbm P^n,\fasm O(d))\lra
H^0(\sbm P^n,\fasm O_Y(d))$ is surjective for $d\geq d^{\prime\prime}$.
\item If (i) is injective for $d^{\prime}=d$ and (ii) is surjective for	
$d^{\prime\prime}=d+1$, then $Y$ has maximal rank.
\end{list}
\end{prop}

\noindent{\bf Proof.}  	See [H1] .$\quad\quad\Box$
\begin{cor} If $Y$ is a zero dimensional closed subscheme of $\sbm P^n$
which is adjusted in degree $d$, then $Y$ has maximal rank.\end{cor}

\noindent{\bf Proof.}  	Put $Y^{\prime}=Y=Y^{\prime\prime}$ in the
proposition.
$\quad\quad\Box$

We will apply 2.2. and 2.3. in the following way. We first define

\begin{defn} For $n,d \geq1 $ let $A_{n,d}$ , $B_{n,d}$ be the integers defined
by 
$$\binm{n+d}{d}= (n+1)A_{n,d} + B_{n,d}$$ 
where $0 \leq B_{n,d}\leq n$. Let $x_i$ ($i = 1,\ldots , A_{n,d}$ ) be generic
closed points of $\sbm P^n$ and let $(y_{n,d} ,L_{n,d} )$ be the 
generic couple
formed by a closed point $y_{n,d}$  and a linear subspace $L_{n,d}$  of
dimension  $B_{n,d} -1$ incident with $y_{n,d}$ . We define $Y_{n,d}$ to be
the union of the first infinitesimal neighbourhoods of
$x_1,\ldots,x_{A_{n,d}}$ in $\sbm P^n$ and the first infinitesimal
neighbourhood of $y_{n,d}$  in $L_{n,d}$ .
\end{defn}
	
We now have

\begin{lem} With the notation of 2.4., $Y_{n,d}$ is numerically adjusted in
degree $d$ and for any $r\geq 0$, there is a uniquely determined $\delta\geq
0$ such that the generic union $Y$ of $r$ first infinitesimal neighbourhoods
of
$\sbm P^n$ is inserted in a sequence of inclusions $Y_{n,\delta} \subset Y
\subset Y_{n,\delta +1}$. If $Y_{n,\delta}$  and $Y_{n,\delta +1}$ are
adjusted then Y has maximal rank.
\end{lem}

\noindent{\bf Proof.} This is a purly numerical consequence of the
definitions and a simple application of 2.2. (c).  $\quad\quad\Box$

We now show that the following proposition implies theorem A.

\begin{prop} With the notation of 2.4, $Y_{n,d}$ is adjusted in degree $d$ in
$\sbm P^n$ for $d\geq 5$. \end{prop}

\noindent{\bf Proof. (of theorem A.)} Let $x_1,\ldots ,x_r$ be closed points
in general position in projective space $\sbm P^n$ and let $Y$ be the union of
the first infinitesimal neighbourhoods of the $x_i$. Let $\delta$ be the
uniquely determined integer such that $A_{n,\delta} \leq r <
A_{n,\delta+1}$. If $d\geq 5$, then $Y$ has maximal rank by 2.5 and 2.6.
If
$\delta<5$, then $Y \subset Y_{n,5}$ so that $Y$ has maximal rank in all
degrees
$d\geq 5$ by 2.2 (b)  and 2.6.
$\quad\quad\Box$
\section{ Preliminaries}
Generally speaking the proof of 2.6 goes as follows (see[H1] for the general
setting). We argue by induction on the dimension $n$ and the degree $d$ using
the following `lemma-definition'

\begin{lemd}. Let $Y$ be a closed subscheme of $\sbm P^n$ numerically
adjusted in degree $d$ and let $H$ be a hyperplane of $\sbm P^n$. We put
$Y^{\prime\prime}= Y\cap H$, which we call the trace of $Y$ on $H$ and we
let
$Y^{\prime}$ be the closed subscheme of $P^n$ defined by the following
canonical exact sequence
\begin{equation}0 \lra I_{Y^{\prime}}(-1) \lra I_Y \lra
I_{Y^{\prime\prime},H}
\lra 0\end{equation}
where $I_Y$ is the ideal of $Y$ in $\sbm P^n$ and $I_{Y^{\prime\prime},H}$ is
the ideal of $Y^{\prime\prime}$  in $H$. We say that $Y^{\prime}$ is the
residual of $Y$ with respect to $H$. In view of definition 1.1, it follows
immediately that $Y$ is  adjusted in $P^n$ in degree $d$ if
the following two conditions are verified
\newcounter{mist1}
\begin{list}{\textbf{ (\roman{mist1}) }}
{\setlength{\leftmargin}{1cm}
\setlength{\rightmargin}{1cm}
\usecounter{mist1}  } 
\item $Y^{\prime\prime}$ is adjusted in $H= \sbm P^{n-1}$ in degree $d$ 
\item	$Y^{\prime}$ is adjusted in $\sbm P^n$ in degree $d-1$. 
\end{list}
\end{lemd}

	 This allows us to resolve the proposition `` $Y$ is adjusted in $\sbm P^n$
in degree $d$ '', into two similar statements, one of which is in lower degree
and the other in lower dimension. Note that to exploit 3.3 in this way, it
is essential to deal only with (numerically) adjusted subschemes!

\subsection{Adjusting} 

Now in practice, $Y$ is the generic member of a nice smooth 
irreducible flat familly $T$. In this case we can use the upper
semicontinuity of the fonctions $t\mapsto h^i(\sbm P^n,\fasm I_{Y_t}(d))$,
$t\in T$, to conclude that
$Y$ is adjusted in
$\sbm P^n$ in degree $d$ if there is some point $t\in T$ where $Y_t$ has
this property (In keeping with common usage we call $Y_t$ a specialisation
of $Y$). This is the essential element which allows us to employ 3.1. Given a
familly $T$ of closed subschemes of $\sbm P^n$ numerically adjusted in
degree $d$, we find a specialisation $Y_t$ and a hyperplane $H$ such that (i)
$Y_{t}^{\prime\prime}= Y_t\cap H$ is numerically adjusted in $H$ in degree
$d$. By the additivity of the Euler characteristic, this implies that (ii) 
$Y_t^{\prime}$ is numerically adjusted in $\sbm P^n$ in degree $d-1$. This
necessity to adjust the Euler characteristic of the trace, is the major
technical obstacle in the proof. Insufficient techniques for doing this was
the source of formidible technical difficulties in the proof of theorem A
given in [A]. The big advance is the differential lemma [AH1] which is given
in a simplified form in 3.3. 
\subsection{ Ejecting}
One technique for adjusting the Euler characteristic 
of the trace is what we call ejecting. Fix a hyperplane $H$ in $\sbm P^n$ and
a point $x$ in some linear subspace $L$ of $H$. Let $X$ be the first
infinitesimal neighbourhood of $x$ in $L$ and let $Z$ be the first
infinitesimal neighbourhood in $\sbm P^n$ of a general point $z$ on a line
$D$ through $x$, with $D$ transverse to $H$. Let $R$ be the unique, closed
subscheme of $\sbm P^n$ obtained as the flat limit of $Z\cup X$, as $z$ is
specialised to $x$ along $D$. The process which consists in replacing
$Z\cup X$ by the specialisation $R$, will be referred to as ejecting  $X$ by
$Z$. The sens of this terminology resides in the fact that the trace of $R$ on
$H$ is the first infinitesimal neighbourhood of $x$ in $H$, while the
residual of R with respect to $H$ is the first infinitesimal neighbourhood of
$x$ in the linear subspace spanned by $L\cup D$ - loosely speaking, in the
limit, $X$ is ejected outside of $H$ and replaced by $Z$.

\subsection{Extending}

\begin{lem} (see [AH2] 1.5) Let $Y,Z$ be a closed subschemes of $\sbm P^n$.
We let
$Y\cap Z$ and $Y\cup Z$ denote respectively the scheme theoretic
intersection and scheme theoretic union of $Y$ and $Z$. If 
$h^i(Z,\fasm I_{Y\cap Z,Z}(d)) =0$ pour $i\geq 0$, then $Y$ is adjusted in
degree
$d$ in $P^n$ if and only if  $Y\cap Z$ is adjusted in $\sbm P^n$ in degree $d$.
\end{lem}

\noindent{\bf Proof.} An immediate consequence of the following, canonical,
exact sequence of ideals
\begin{equation}	0 \lra \fasm I_{Y\cup Z} \lra \fasm I_Y \lra \fasm
I_{Y\cap Z,Z}
\lra 0\end{equation}

The process, in 2.4.1, which consists in replacing $Y$ by $Y\cup Z$,
will be referred to as extending  $Y$ by $Z$.

\subsection{ The differentielle method}

	In \S4, we will use the following lemma [AH1, 1.3, 1.5].

\begin{lem} We consider a closed subscheme of $\sbm P^n$, numerically
adjusted in degree $d$ which is the disjoint union of two closed 
subschemes $W, Z$, with $Z$
of finite support. We suppose that the residual $Z^{\prime}$ of $Z$ with
respect to a fixed hyperplane $H$ is contained in the trace
$Z^{\prime\prime}$ of
$Z$ on
$H$. Then the union $Y$ of $W$ with the generic translate of $Z$ in
$\sbm P^n$ is adjusted 
 in $\sbm P^n$ in degree $d$ if the following two conditions are verified
$$\br{ll}
\mbox{(dime)}&	\mbox{There is a subscheme $d^{\prime\prime}Z$ of
$Z^{\prime\prime}$ containing
$Z^{\prime}$ such that $W^{\prime\prime}\cup d^{\prime\prime}Z$ 		 is
adjusted }\\
&\mbox{in $H$, in degree $d$.}\\
&\\
\mbox{(degue)}&	\mbox{There is a closed subscheme $d^{\prime}Z$
contained in the intersection of the base of the } \\
&\mbox{linear system $|H^0(H,\fasm I_{W^{\prime\prime}\cup
Z^{\prime},H}(d))|$ with
$Z^{\prime\prime}$ ,	such that $W^{\prime}\cup d^{\prime}Z$ is adjusted
in}\\ &\mbox{$\sbm P^n$, in degree
$d-1$.}
\er$$
\end{lem}
\begin{rems}{\em
\be
\item The words dime and degue indicate that the corresponding conditions
are respectively in lower dimension and lower degree.
\item When $Z$ is the first infinitesimal neighbourhood in $\sbm P^n$ of a
closed point in $H$ with
$$\gchi(H,\fasm I_{W^{\prime\prime}\cup Z^{\prime},H}(d))>0\quad ,
\quad \gchi(H,\fasm I_{W^{\prime\prime}\cup Z^{\prime\prime},H}(d))<0$$
it is easy to see that, if
$d^{\prime\prime}Z$ exists verifying dime, then the intersection $\fasm
B$ of the base of the linear system $|H^0(H,\fasm I_{W^{\prime\prime}\cap
Z^{\prime},H} \otimes \fasm L_H)|$ with $Z\upp$, has degree
$\mbox{deg}(Z)-\mbox{deg}(d\upp Z)$. Now this is just the condition that
$W\up \cup \fasm B$ be numerically adjusted in $\sbm P^n$ in degree $d-1$,
so that
$\fasm B=d\up Z$ is the unique choice for $d\up Z$ with $W\up \cup d\up Z$
numerically adjusted in $\sbm P^n$ in degree $d-1$. What is more, this 
$d\up Z$ depends only on $W\upp$ and the point $Z\up$.
	That $\fasm B$ have degree $\mbox{deg}(Z)-\mbox{deg}(d\upp Z)$, just
says that locally at the point $Z\up$, the base of the linear system
$|H^0(H,\fasm I_{W\upp \cup Z\up ,H}\otimes \fasm L_H)|$ is smooth of
codimension $\mbox{deg}(d\upp Z)$ in $H$. As such $\fasm B$ is just the
intersection of $Z\upp$ with some linear subspace $L$ through $Z\up$, with
$L$ transverse to $d\upp Z$ and of dimension $n-deg(d\upp Z)$. In particular,
$\fasm B$ is the first infinitesimal neighbourhood of $Z\up$ in $L$.  
\ee}
\end{rems}
\section{Proof of 2.6}

\subsection{proof of 2.6}

For $d\geq 5$, $Y_{n,d}$ is adjusted in $\sbm P^n$
in degree $d$. 

\noindent{\bf Proof.} We proceed by induction on $n$, the case $n=1$ being
classical one variable interpolation.
	
	Let $H$ be the generic hyperplane of $\sbm P^n$ ($n\geq 2$). There are two
cases to consider in accordance with $B_{n-1,d} = 0$ and $B_{n-1,d} > 0$. 
	
Case $B_{n-1,d} = 0$. 

In this case we let $S_{n,d}$ be the specialisation of $Y_{n,d}$
obtained by specialising the points $x_i$ ($i= 1,...,A_{n-1,d}$) into $H$. Then
the trace $S_{n,d}\upp$ of $S_{n,d}$ on $H$ is just $Y_{n-1,d}$. By the
induction hypothesis $Y_{n-1,d}$ is adjusted in degree $d$, so by 3.1 it is
enough to show that the residual $S_{n,d}\up$ of $S_{n,d}$ with respect to
$H$, is adjusted in degree $d-1$.

Case $B_{n-1,d} > 0$.

In this case we let Sn,d be the specialisation of $Y_{n,d}$, 
 obtained by specialising the points $x_i$ ($i= 1,...,A_{n-1,d +1}$) into $H$.
We apply 2.5.1 with $Z$ the first infinitesimal neighbourhood of $x_1$ and
$d\upp Z$ the intersection of $Z$ with the generic linear subspace of $H$
passing through $x_1$ of dimension $B_{n-1,d -1}$. The subscheme ($S_{n,d}
-Z)\upp \cup d\upp Z$ is just $Y_{n-1,d}$ so that the dime of 2.5.1 results
from the induction hypothesis. As indicated in 2.5.2 (2), there is a unique
choice for $d\up Z$ verifying the numerical part of the condition degue of
4.1, we must show that the union of $S_{n,d} -Z)\up$ and $d\up Z$ is
adjusted in degree $d-1$.

	Now we define $G_{n,d}$ to be $S_{n,d}\up$ if $B_{n-1,d} = 0$ (resp.
$(S_{n,d} -Z)\up$ if $B_{n-1,d} > 0$). This is just the generic union of
$A_{n,d} - A_{n-1,d}$  (resp. $A_{n,d} - A_{n-1,d -1}$) double points in
$\sbm P^n$ with $A_{n-1,d}$ points in $H$. To complete the two cases 
it will be enough to prove that if $Z$ is a generic double point of $H$ and
$Z_0$ is any subscheme of $Z$ such that $T_{n,d} = G_{n,d}\cup Z_0$ is
numerically adjusted in degree $d-1$, then $T_{n,d}$ is adjusted in degree
$d-1$. This will be proven in ¤4. This result is more general than the degue,
but, as we do not have any useful information about $d\up Z$, we are obliged,
in practice, to prove the more general result. $\quad\quad\Box$
\section{The Degue }
In the proof of the following proposition we simply show how the
proposition can be deduced from results in lower degree. These results in
lower degree are then proven in ¤6.

\begin{prop} For $d\geq 5$, let $G_{n,d}$ be the generic union of
$A_{n,d} - A_{n-1,d} - \varepsilon$ double points of $\sbm P^n$ with
$A_{n-1,d}$ points in $H$, where $\varepsilon =0 $ if $B_{n-1,d} = 0$ and
$\varepsilon =1$ if $Bn-1,d > 0$. Let $Z$ be a generic double point of $H$ and
let $Z_0$ be any subscheme of $Z$ such that $T_{n,d} = G_{n,d} \cup Z_0$ is
numerically adjusted in degree $ d-1$, then $T_{n,d}$ is adjusted in degree
$d-1$.
\end{prop}

\noindent{\bf Proof.} If $n=1$ this follows by classical one variable
interpolation. We will now treat the cases $n=2$ and $n\geq2$ separately.

\noindent Case (1). n=2.

In this case we have $B_{2,d} = 0$ or 1, since $(d+2)(d+1)/2 \neq 0$ or 1
(mod 2). 
If $d=5$, then $B_{1,5} = 0 = B_{2,5}$ and $T_{2,5}$ is adjusted by 5.1. 
If $d=6$ then $T_{2,6}$ is the union of five first infinetisimal
neighbourhoods and one closed point in $\sbm P^2$ with three closed points
and one infinetisimal neighbourhood in $H= \sbm P^1$. By ejecting one
closed point in $H$ with one infinitesimal neighbourhood of $\sbm P^2$, the
resulting specialisation $S$ has an adjusted trace on $H=\sbm P^1$. The
residual
$S\up$ of $S$ with respect to $H$, can be specialised to the generic union of
four infinitesimal neighbourhoods of $\sbm P^2$ with one closed point in $H$
and one degree two closed subscheme of $\sbm P^2$, transverse to $H$ with
support a closed point of
$H$. It follows from 5.1 that $S\up$ is adjusted in $\sbm P^2$ in degree
$d=4$. 

\begin{lemb} Let $s,t$ be the unique integers satisfying 
$$2s-t = d - (d+1)/2 + B_{1,d} /2 - B_{2,d} - B_{1,d}(3- B_{1,d}) \quad
;\quad   0
\leq t
\leq1.$$ 
Then for $d\geq 7$, we have $(s-t)\geq 1$ et $s\leq d/2 +1/2$.
\end{lemb}

\noindent{\bf Proof. (of lemma)} If $d=7$ then $B_{1,d} = 0 = B_{2,d}$ and
$s= 2, t=1$. If $d=8$, we have $B_{1,d} = 1$, $B_{2,d} = 0$ giving $s=1$,
$t=0$. So we have the result in these cases. For $d\geq 9$, we consider the
expression

\begin{equation}	2(s-t) = d/2 - 1/2 + B_{1,d} /2 - B_{1,d}(3- B_{1,d}) -
B_{2,d }- t\end{equation}

Noting that $B_{1,d}, B_{2,d}$ and t are all less 
than or equal to one, we find
$$2(s-t) \geq 1/2 > 0$$
which gives the first part of the lemma. The second part is trivial.
$\quad\quad\Box$

	Now suppose that $d\geq 7$ and let $s,t$ be as in lemma 5.1.1. Let $S$ be
the specialisation of $T_{2,d}$ obtained by specialising $(s-t)\geq 0$ first
infinitesimal neighbourhoods to have support in $H$, ejecting $t\leq1$ of the
$A_{1,d} \geq 1$ closed points of $H$ with a further $t$ infinitesimal
neighbourhoods and finally specialising $y_{2,d}$ to $H$ in the case $B_{2,d}
=1$. By (5.1.1), the trace of $S$ on $H$ is numerically adjusted in degree $d$
in $H= \sbm P^1$, hence it is adjusted. Now the residual $S\up$ of $S$ with
respect to $H$ is the numerically adjusted, generic union of a number of
infinitesimal neighbourhoods of $\sbm P^2$ with $(s-t)\geq1$ closed points
in $H$ and $t$ finite degree two closed subschemes transverse to $H$ with
support a closed point of $H$. Since, by 5.1.1, $s\leq d/2+1/2$ , we
can conclude by 6.1, that $S\up$ is adjusted in $\sbm P^2$ in degree $d-2$.
This gives the result by 3.1.  

\noindent Case (2)  $n\geq 3$.

	The scheme $T_{n,d}$ contains the subscheme $d\up Z$ 
 which depends on the $A_{n-1,d}$ closed points in the trace  $T_{n,d}\up$
(see 3.4 (b)). With $n\geq 3$ we have very little information about the linear
subspace $L$ spanned by $d^{\prime}Z$ other than the fact that it is
transverse to $d\upp Z$. Any specialisation of $T_{n,d}$ which modifies the
closed points in $H$ will modify the subscheme $d\up Z$. We avoid this
unknown quanity by ejecting $d\up Z$ by one of the first infinitesimal
neighbourhoods in $T_{n,d}$. We then specialise further first infinitesimal
neighbourhoods to
$H$, ejecting a number of closed points, in such a way that the trace of this
specialisation is numerically adjusted in $H$ in degree $d-1$. In this way
$d'Z$ does not appear in the trace and, when it appears in the residual, it no
longer depends on the closed points in $H$.

	As usual we begin with an accounting lemma.

\begin{lemb} For $n\geq 3$, $d\geq 5$, let $s,t$ be the integers determined 
by the expression

\begin{equation}ns-t = \binm{n+d-2}{d-1} - A_{n-1,d} - B_{n,d} -
n\delta_{n-1,d}
\quad ;\quad 0\leq t<n\end{equation}

where $\delta_{n-1,d} =0$ if $B_{n-1,d} = 0$ and  $\delta_{n-1,d} =1$
otherwise. Then we have the following inequalities:
\newcounter{mist2}
\begin{list}{\textbf{ (\alph{mist2}) }}
{\setlength{\leftmargin}{1cm}
\setlength{\rightmargin}{1cm}
\usecounter{mist2}  } 
\item	$s-t \geq (n-1) + \delta_{n-1,d} \,B_{n-1,d }-1$
\item	$s    \geq n(n-1)/2 + \delta_{n-1,d} \,B_{n-1,d} -1$

\noindent or one of the following
$$\br{l}	n=4,\;  d=5 \mbox{ where } s=5,\;  t=0,\;  B_{4,5} =1,\; B_{3,5} =
0\\
n=5, \; d=5 \mbox{ where } s=8, \; t=0, \; B_{5,5} = 0 \\
n=6, \; d=5 \mbox{ where } s=14, \; t=0,\;  B_{6,5} = 0
\er$$
\item $\frac{1}{n}\binm{n+d-2}{d-1} + (n-1)/2 \geq s$
\end{list}
Equally,
$$A_{n-1,d} \geq n(n-1)/2 + d_{n,d} (n- B_{n,d})$$
  
\end{lemb}

\noindent{\bf Proof. ( of the lemma)} We have 
$A_{n-1,d} = \frac{1}{n}\binm{n+d-1}{d} - B_{n-1,d} /n$, so that
\begin{equation}	ns-t = \frac{(n-1)(d-1)}{nd} \binm{n+d-2}{d-1}
- B_{n,d} - n\delta_{n-1,d}+ \frac{B_{n-1,d}}{n}\end{equation}
Since $t\leq (n-1)$, we see that (a) is implied by (b) for $n\geq 4$. 
We will prove (a) for $n=3$. From (5) we get 
$$3(s-t)=\frac{(d+1)(d-1)}{3} - B_{3,d} -3 \delta_{2,d} + \frac{B_{2,d}}{3} -
2t$$
so that (a) is equivalent, in its strict inequality form, to 
\begin{equation} \br{lcl}\frac{(d+1)(d-1)}{3} & > &3(1+ \delta_{2,d} 
B_{2,d-1}) + B_{3,d} + 3 \delta_{2,d} - \frac{B_{2,d}}{3} + 2t\vspace{.5ex}\\
				     & = &3 + B_{3,d} + (3 \delta_{2,d}- 1/3) B_{2,d} +
2t\er\end{equation}
Using the inequalities $t \leq 2$, $B_{3,d} \leq 3$,
$B_{2,d} \leq 2$ we see that it is enough to show that
	     $$(d+1)(d-1)/3  > 46/3$$
This is the case for $d\geq 7$. For $d=5$ we have $\delta_{2,5}=B_{2,5}
=B_{3,5} =0$ and $s=3$, $t=1$. For $d=6$ we have $\delta_{2,6}=B_{2,6} =1$,
$B_{3,6} =0$ and $s=3$, $t=0$. This gives (a) for $n=3$ as well.

We will now prove (b). By (5) we have
\begin{equation}	ns=
\frac{n-1)(d-1)}{nd}\binm{n+d-2}{d-1} 
					- B_{n,d} - n\delta_{n-1,d} + \frac{B_{n-1,d}}{n} + t
\end{equation}
so that (b) is equivalent to
$$\br{lcl}
	\frac{n-1)(d-1)}{nd}\binm{n+d-2}{d-1}
	&> &\frac{n^2(n-1)}{2} + n \delta_{n-1,d} (B_{n-1,d}  - 1) - n + B_{n,d}
-\frac{B_{n-1,d}}{n} + n.\delta_{n-1,d} - t \vspace{.5ex}\\
	&= &\frac{n^2(n-1)}{2} + B_{n,d} + \left(n\delta_{n-1,d}-\frac{1}{n}\right)
B_{n-1,d}- n -t
\er$$
Using the upper bounds for $B_{n,d}$ , $B_{n-1,d}$ one easily sees that it is
sufficient to prove the following inequality
	$$	\frac{(n-1)(d-1)}{nd}\binm{n+d-2}{d-1} 
		> \frac{n^2(n-1)}{2} + \frac{(n+1)(n-1)^2}{n}$$
Simplifying, this becomes
\begin{equation}(d-1)\binm{n+d-2}{d}>
\frac{1}{2}(n-1)(n^3+2n^2-2)\end{equation}
and one easily verifies (8) for $d=7$, $n\geq 3$; $d=6$, $n\geq 6$ and
$d=5$,
$n\geq11$. Now using the fact that the left hand side of (8) is an
increasing function of $d$, we conclude that to prove (b) it is enough to treat
the particular cases $d=6$, $n=3,4$; $d=5$, $n=3,\ldots ,10$.

For $d=6$ we have we have for $n=3$ (resp. $n=4$) $s=3$, $t=0$ (resp.
$s=9$, $t=1$) giving (b) in this case.
	For $d=5$ we have the following triples $(n,s,t)$ : $(3,3,1)$, $(4,5,0)$,
$(5,8,0)$, $(6,14,0)$, $(7,21,4)$, $(8,29,1)$, $(9,39,1)$, $(10,51,5)$, giving
(b) in these cases as well.
	
One easily proves (c) using the expression (7). 

Finally to prove the last part of the lemma we have 
$$A_{n-1,d} =\frac{1}{n}\binm{n+d-1}{d} - B_{n-1,d}/n $$
so we must show that
$$\frac{1}{n}\binm{n+d-1}{d} > \frac{n(n-1)}{2} + \delta_{n,d} (n- B_{n,d}) +
\frac{B_{n-1,d}}{n} -1$$
Using the upper bound for $B_{n-1,d}$ and the lower bound  $1\leq B_{n,d}$
for
$\delta_{n,d}\neq 0$, we see that, since the left hand side is an increasing
function of $n$, it is enough to show that
$$\binm{n+4}{5} > \frac{n^2(n-1)}{2} + n(n-1)$$
or that
$$(n+4)(n+3)(n+1) > 60(n-1)$$
One easily verifies that this holds for $n\geq 3$.$\quad\quad\Box$

Now let $S$ be the specialisation of $T_{n,d}$ obtained by specialising
$(s-t)$ first infinitesimal neighbourhoods to have support in $H$,
specialising $(y_{n,d} ,L_{n,d})$ into $H$, then ejecting $t\leq n-1$ closed
points of H by a further t infinitesimal neighbourhoods (one easily verifies
that one has enough). In the case where $d\up Z$ is non empty we eject
$d\up Z$ by yet another infinitesimal neighbourhood. By (4) the trace
$S\upp$ of $S$ on $H$ is numerically adjusted in degree $d-1$.  We will
apply lemma 3.1 to $S$.

We first show that $S\upp$ is adjusted in degree $d-1$ in $H= \sbm
P^{n-1}$. For this it is enough to show that $S\upp$ admits a specialisation
to a closed subscheme of  $H$ of the type figuring in 6.1.

If $B_{n,d}=0$, then $S\upp$ is a generic union of closed points and first
infinitesimal neighbourhoods of closed points of $H$. If $B_{n,d}\neq 0$, then
we can specialise $\delta_{n,d}(n- B_{n,d})$ of the closed points in $H$ to
the point $y_{n,d}$ , completing the first infinitesimal neighbourhood of
$y_{n,d}$ in $L_{n,d}$ to a first infinitesimal neighbourhood of $y_{n,d}$ in
$H$. This means that $S\upp$ admits a specialisation which is a generic
union of first infinitesimal neighbourhoods of closed points in H and a
number $c$ of closed points of $H$ where 
$$c = A_{n-1,d} -  \delta_{n,d}(n- B_{n,d}) -t \geq \frac{(n-1)(n-2)}{2}$$
By specialising further closed points to infinitesimal neighbourhoods we
obtain a specialisation where the number $c$ of closed points satisfies
$$\frac{(n-1)(n-2)}{2} + (n-1) \geq c \geq \frac{(n-1)(n-2)}{2}$$
and we conclude that $S\upp$ admits a specialisation to a closed
subscheme of the type figuring in 6.1, so that $S\upp$ is adjusted in degree
$d-1$.
	
We will now show that $S\up$ is adjusted in degree $d-2$ in $\sbm P^n$. The
residual  $S\up$ of $S$ with respect to $H$ is the generic union, numerically
adjusted in $\sbm P^n$ in degree $d-2$, of first infinitesimal
neighbourhoods of closed points in $\sbm P^n$ with $t\leq n-1$ finite,
degree two, closed subschemes with support in $H$, but transverse to $H$,
$s-t$ closed points of $H$ and $\mbox{res}(d\up Z)\subset H$ (this is the
residual of the ejected $d\up Z$ which we consider as empty if $B_{n-1,d}$
is zero). 
	
Now the subscheme $\mbox{res}(d\up Z)$ is independent of the 
$$s-t \geq (n-1) + \delta_{n-1,d} (B_{n-1,d} - 1) $$
closed points in $H$, so if $\mbox{res}(d\up Z) $ is non-empty we can
specialise $(B_{n-1,d}  -1)$ of these closed points to $d\up Z$ to make a
first infinitesimal neighbourhood. The conditions of lemma 5.1.2, assure
that this specialisation of $S\up$ satisfies the hypotheses of 6.1, where the
variable $v$ is zero if $B_{n-1,d}$ is zero and one otherwise. As such $S\up$
is adjusted in $\sbm P^n$ in degree $d-2$ and we have completed the proof
of 5.1. $\quad\quad\Box$
\section{Results and initial cases}
\begin{prop}  Let $H=\sbm P^{n-1}$  be a hyperplane of $\sbm P^n$. In
$\sbm P^n$ we consider, for $d\geq 3$, the closed subscheme $Y$ formed by
the generic union, numerically adjusted in $\sbm P^n$ in degree $d$, of 
\newcounter{mist3}
\begin{list}{\textbf{ (\alph{mist3}) }}
{\setlength{\leftmargin}{1cm}
\setlength{\rightmargin}{1cm}
\usecounter{mist3}  } 
\item the first infinitesimal neighbourhoods of $u$ closed points of
$\sbm P^n$
\item	$a$ closed points of H,
\item	$b$  finite closed subschemes of degree two of $\sbm P^n$ transverse
to $H$ with support a closed point of $H$,
\item	the first infinitesimal neighbourhood in $\sbm P^n$ of $v\leq 1$
closed points of
$H$ where
\newcounter{mist4}
\begin{list}{\textbf{ (\roman{mist4}) }}
{\setlength{\leftmargin}{1cm}
\setlength{\rightmargin}{1cm}
\usecounter{mist4}  } 
\item	$a \geq n-1$, or $d=3$ and we have one of the following cases
$$n=4, a=5, b=0=v \quad ;\quad
		n=5, a=8, b=0, v=1\quad ;\quad 
		n=6, a=14, b=0=v $$
\item $\frac{1}{n}\binm{n+d}{d+1}+ \frac{n-1}{2} \geq a+b
\geq \frac{n(n-1)}{2}$
\item $0 \leq b < n$
\item if $n=2$ and $d=3$, then $v=0$
\end{list}\end{list}
Then $Y$ is adjusted in $\sbm P^n$ in degree $d$.
\end{prop}

\noindent{\bf Proof. } The proposition is true for $n=1$ by classical one
variable interpolation. For $d=3$ the proposition results from 6.3 below.
Henceforth we suppose that $n\geq 2$ and $d\geq 4$ and we argue by
induction on $n$ and $d$. Firstly we need the following accounting lemma.

\begin{lemb} With the same hypotheses as in 6.1, let s,t be the integers
uniquely determined by the following conditions,
\begin{equation}sn - t = (n+d-1)!/(n-1)! d! - (a+b)   ;   	0 \leq t < n
\end{equation}
Then for $n\geq 2$, we have
\newcounter{mist5}
\begin{list}{\textbf{ (\alph{mist5}) }}
{\setlength{\leftmargin}{1cm}
\setlength{\rightmargin}{1cm}
\usecounter{mist5}  } 
\item	$s-t+b \geq n-1$
\item	$s+b    \geq \frac{n(n-1)}{2}$, or $n=4, a=15, b=0$ in which case
$t=0$ and $s=5$.
\end{list}
\end{lemb}

\noindent{Proof. (of lemma) } For (a) we have
$$s-t = \frac{1}{n}\binm{n+d-1}{d} - \frac{(a+b)}{n} -\frac{(n-1)t}{n}$$
and using (ii) of the proposition,we have
\begin{equation}\br{lcl}	s-t &\geq& \frac{1}{n}\binm{n+d-1}{d} - 
\frac{1}{n^2}\binm{n+d}{d+1} -
\frac{n-1}{2n} -\frac{(n-1)t}{n}\\
&=& \frac{(n-1)d}{n^2(d+1)}\binm{n+d-1}{d} - \frac{n-1}{2n} 			
-\frac{(n-1)t}{n}\\
& =& R(n,d,t)
\er\end{equation}
We want $(s-t) > (n-2)$ . Since $t\leq (n-1)$ and $R(n,d,t)$ is a strictly
increasing function of $d$, it is enough to show that $R(n,4,n-1) > (n-2)$, i.e. 
$$\frac{4(n-1)}{5n^2}\binm{n+3}{4} > \frac{4n-7}{2}$$
This is the case for $n=2$ or $n\geq 4$. 

	If $n=3$, then (10) becomes $s-t \geq \frac{7-2t}{3} \geq 1$ and the
last inequality is strict for $t<2$. When $t=2$, $s-t\geq $1 with strict
inequality if $s\geq 3$. Finally if $s=3$ and $t=2$, we find by (9) that
$(a+b) = 8$. Using the fact that $Y$ is adjusted in $\sbm P^3$ in degree 
three, we find $b+8+4(u+v)=35$, so that $b\geq 3$ contradicting (iii) of the
proposition. This gives (a).

For (b) we have 
\begin{equation}s = \frac{1}{n}\binm{n+d-1}{d} -
\frac{a+b}{n} +
\frac{t}{n}
\end{equation}
and using (ii) of the proposition we find
\begin{equation}\br{lcl}	s &\geq& \frac{(n-1)d}{n^2(d+1)}\binm{n+d-1}{d}
- \frac{n-1}{2n} + \frac{t}{n}\\
&\geq& \frac{4(n-1)}{5n^2}\binm{n+3}{4}- \frac{n-1}{2n}
+ \frac{t}{n}
\er
\end{equation}
Since $t\geq 0$ and we need $s+b > \frac{n(n-1}{2}-1$, it is enough to show
that 
$$\frac{4(n-1)}{5n^2}\binm{n+3}{4} > \frac{n(n-1)}{2}-1- \frac{n-1}{2n}$$
This is the case for $n=2,3$ or $n\geq 5$. 
	
For $n=4$, (12) becomes $s \geq \frac{39 + 2t}{8}$, so that $s>5$ for
$t>0$. If $t=0$, the same expression gives $s\geq 5$. Finally if $s=5$ and
$t=0$, then (11) gives $a+b = 15$. If $b>0$, then $s+b\geq 6$ as required
and we are left with the case $a=15$, $b=0=t$, $s=5$. This gives (b) using
the fact that the right hand side of (12) is a strictly increasing function
of d. $\quad\quad\Box$

Now to finish the proof of the proposition we use lemma 2.1. 
Let $s,t$ be as in the (9). We consider the specialisation $S$ of $Y$
obtained by specialising $(s-t)-v$ ($\geq 0$ by the lemma) of the $u$ first
infinitesimal neighbourhoods to have support in $H$ and ejecting $t$ of the
$a\geq (n-1)$ points in $H$ by a further $t$ of the infinitesimal
neighbourhoods. By (9), the trace $S\upp$ of $S$ on $H$ is numerically
adjusted in $H$ in degree $d$ with 
\begin{equation}a\up = a+b-t \geq \frac{(n-2)(n-1)}{2} \end{equation}
closed points contained in $S\upp$. 

For $n\geq 2$, $S"\subset H=\sbm P^{n-1}$ can be specialised to a 
subscheme of the same type as figures in the proposition. In fact, $S\upp$ is
a generic union of closed points and first infinitesimal neighbourhoods of
$\sbm P^{n-1}$ and by specialising closed points to infinitesimal
neighbourhoods, we can suppose that the number of closed points satisfies
\begin{equation} \frac{(n-2)(n-1)}{2} + (n-1)  \geq   a\up 
\geq  \frac{(n-2)(n-1)}{2}\end{equation}
One then verifies that 
\begin{equation} \binm{n+d-2}{d}+ \frac{n-1}{2} \geq \frac{n(n-1)}{2 }+
\frac{n-1}{2}
						           \geq \frac{(n-2)(n-1)}{2} + (n-1)\end{equation}
for $d\geq 2$. Now we conclude by the induction hypothesis that $S\upp$ is
adjusted in degree $d$ in $\sbm P^{n-1}$.

The residual $S'$ of $S$ with respect to $H$ is numerically adjusted in
$\sbm P^n$ in degree $d-1$. It is the generic union of $(s-t+b)$ closed points
in $H$, plus $t\leq (n-1)$ finite, degree two, closed subschemes of
$\sbm P^n$ with support in $H$, but transverse to $H$ (these are the
residuals of the ejections) and $(u-s)$ first infinitesimal neighbourhoods of
closed points in $\sbm P^n$. In view of (9), (ii) and (iii) of the
proposition and lemma .1.1, we have
$$(s-t+b) \geq (n-1) 
		\frac{1}{n}\binm{n+d-1}{d} + \frac{n-1}{2} \geq (s+b) \geq (n-1) $$
so that $S\up$ is of the type figuring in the theorem for $d\geq 4$. By the
induction hypothesis we conclude that $S\up$ is adjusted in degree $d-1$ for
$d\geq4$. 

The proposition now follows by 3.1. $\quad\quad\Box$
\begin{defn}  Let $x_1, \ldots ,x_{s+2}$ be closed points in $\sbm P^n$ 
with $s\leq n-1$ and consider the union $U$ of the $\frac{(s+2)(s+1)}{2}$
lines which join these points two by two. If $H$ is a hyperplane in
$\sbm P^n$ not containing any of the points, then we say that the
intersection $U\cap H$ is an $s$-complex. The points of an $s$-complex span
a linear subspace of dimension $s$ and we will admit the empty set as a
(-1)-complex.
\end{defn}
\begin{prop} Let $Y$ be the generic union, numerically adjusted in $\sbm
P^n$ in degree $d=3$, of
\newcounter{mist6}
\begin{list}{\textbf{ (\alph{mist6}) }}
{\setlength{\leftmargin}{1cm}
\setlength{\rightmargin}{1cm}
\usecounter{mist6}  } 
\item the first infinitesimal neighbourhoods of $u$ closed points of $\sbm
P^n$
\item an $(l-2)$-complex 
\item $a$ closed points of $H$
\item $b$  finite closed subschemes of degree two of $Pn$ transverse to 		
 $H$ with support $a$ closed point of $H$,
\item the first infinitesimal neighbourhood in $\sbm P^n$ of $v\leq 1$ 
 closed points of $H$, with $v=0$ If $n=2$, where
\newcounter{mist7}
\begin{list}{\textbf{ (\alph{mist7}) }}
{\setlength{\leftmargin}{1cm}
\setlength{\rightmargin}{1cm}
\usecounter{mist7}  } 
\item $a \geq n-1$
\item $\binm{n+3}{4} + \frac{n-1}{2} \geq a+b \geq \frac{n(n-1)}{2}$ 
or one of the following cases
$$n=4, a=5, b=0=v, l=-1$$
$$n=5, a=8, b=0, v=1, l=-1$$
$$n=6, a=14, b=0=v, l=-1$$
\item $0 \leq b < n$
\item $-1 \leq l \leq \frac{n+2}{2}$
\end{list}\end{list}
Then $Y$ is adjusted in $\sbm P^n$ in degree $d=3$.
\end{prop}

\noindent{\bf Proof.}  For $n=1$ this is just classical one variable
interpolation, so from now on we suppose that $n\geq 2$ and we argue by
induction on $n$ using lemma 3.1.
	We begin with the following accounting lemma.

\begin{lemb}  With the hypotheses of the proposition, let $w=u+v$,
then there exists two integers $s,t$ verifying the following conditions
\newcounter{mist8}
\begin{list}{\textbf{ (\alph{mist8}) }}
{\setlength{\leftmargin}{1cm}
\setlength{\rightmargin}{1cm}
\usecounter{mist8}  } 
\item[$(\alpha)\quad$] $(w-s)n - t = \binm{n+2}{3} - (a+b) - \binm{s}{2}$
\item[$(\beta)\quad$]	$0 \leq t < n$
\item[$(\gamma)\quad$] $1 \leq s \leq \frac{n+1}{2}$
\item[$(\delta)\quad$] $(w-s) \geq (t+v)$
\end{list}
\end{lemb}

\noindent{\bf Proof. (of lemma)} We begin with the exceptional cases of
6.3 (ii). For $n=4$, $a=5$, $b=0=v$, $l=-1$ we have $u=w=6$ and it is enough
to take $s=2$ and  $t=2$. For $n=5$, $a=8$, $b=0$, $v=1$, $l=-1$ we have
$u=7$ and it is enough to take $s=3$, $t=1$. For $n=6$, $a=14$, $b=0=v$,
$l=-1$ we have $u=w=10$ and it is enough to take $s=3$, $t=3$.
	Now, in the general case, we consider the following function of $s$,
\begin{equation}	
F(s) = \frac{(n+2)(n+1)n}{6} - (w-s)n - (a+b) - \frac{s(s-1)}{2}
\end{equation}
which is a decreasing fonction of s on the intervalle $[1, ... ,n]$, with
$0 \leq F(s+1) - F(s) = n-m \leq n-1$.

As such, there is an s verifying $(\alpha)$, $(\beta)$ and $(\gamma)$ if we
know that the following two conditions are satisfied,
$$		(1) \quad F(1) \leq 0         \quad\quad         (2) \quad 
F[\![\frac{n+1}{2}]\!]
\geq -n+1 $$
where $[\![\;\; ]\!]$ denotes the integer part.
The fact that $Y$ is adjusted in degree $d=3$ in $\sbm P^n$, allows us to
express $w$ as
$$w = (n+3)(n+2)/6 - (a+b)/(n+1) -l(l-1)/2(n+1)$$
and this gives
$$	F(s) = -n(n+2)/3 - (a+b)/(n+1) + nb/(n+1) + 
ns + nl(l-1)/2(n+1) - s(s-1)/2
$$
Now we have
$$(n+1)(F(1)-1) =  -n(n+1)(n+2)/3 - (a+b) + nb + n(n+1) + nl(l-1)/2 - (n+1)$$
and using (ii) and (iv) of the proposition we find
$$\br{lcl}(n+1)(F(1)-1) &	\leq &-n(n+1)(n+2)/3 - n(n-1)/2 + n(n-1) 
+ n(n+1)+\vspace{.5ex} \\
&&+
n^2(n+2)/8 - (n+1)\vspace{.5ex}\\ 
&=& (-5n^3 + 18n^2 - 28n - 24)/24 \\
&<& 0 
\er$$
pour $n\geq2$.
On the other hand,
$$\br{lcl}(n+1)(F(n/2)+n) &= &-n(n+1)(n+2)/3 - (a+b) + bn+ n^2(n+1)/2
\vspace{.5ex}\\&&+ nl(l-1)/2- n(n+1)(n-2)/8 + n(n+1)\vspace{.5ex}\\
&\geq &-n(n+1)(n+2)/3 -((n+3)(n+2)(n+1)/24 +\vspace{.5ex}\\+ (n-1)) +
           n^2(n+1)/2 \vspace{.5ex}\\&&- n(n+1)(n-2)/8 + n(n+1)\vspace{.5ex}\\
&= &(9n^2 - 9n + 6)/24 \\
&>& 0 
\er$$
pour $n\geq 1$.

Finally we must show that if $s$ satifies $(\alpha)$, $(\beta)$ and
$(\gamma)$, then
$(w-s)\geq (t+v)$. By $(\alpha)$, we have, using (ii) and (iv) of the
proposition, 
$$\br{lcl}w-s &=& \frac{1}{6}(n+2)(n+1) - \frac{1}{n}(a+b) -
\frac{1}{2n}s(s-1) + \frac{t}{n}\vspace{.5ex}\\ &\geq &\frac{
1}{6}(n+2)(n+1) -
\frac{1}{n}\left(\frac{1}{24}(n+3)(n+2)(n+1) + \frac{1}{2}(n-1)\right) -
\frac{1}{8n}(n^2-1) +
\frac{t}{n}\vspace{.5ex}\\ &= &\frac{3}{4}(n+2)(n+1)(n-1) - \frac{n}{8} +
\frac{1}{2n} +
\frac{t}{n}\vspace{.5ex}\\ &=& V(n) + \frac{t}{n}
\er$$
We will show that $V(n) + t/n > t+v-1$. Since $(n-1)\geq t$, it will 
be enough to show that 
$$W(n) =  V(n) - ((n-1)2/n) -1 > v$$

For $n=2$, we have $v=0$ and $W(2)=1$. Now 
$W(n) = \frac{1}{8}(n^3 + 2n^2 - 10n + 14)  - \frac{1}{2n}$, is a strictly
increasing function of $n$ for $n\geq2$ so that, since $v\leq 1$, $W(n)>v$
for $n\geq 2$. This completes the proof of the lemma. $\quad\quad\Box$
	
Now to complete the proof of the proposition we use lemma 2.1. 
Fix $s,t$ verifying the conditions of lemma 6.3.4. Let $S$ be the
specialisation of $Y$ obtained by specialising $u-s-t$ (which is $\geq 0$ by
lemma 6.3.4 (d)) of the $u$ first infinitesimal neighbourhoods to have
support in $H$, and ejecting $t$ of the a closed points in $H$ with a further
$t$ infinitesimal neighbourhoods. We will show that $S$ is adjusted in
$\sbm P^n$ in degree $d=3$.

Let $S_0$ be obtained from $S$ by extending the $s$ infinitesimal
neighbourhoods remaining outside of $H$ by the union $Z$ of the $s(s-1)/2$
lines which join them two by two. Clearly $S\subset D$ is adjusted in each
of the lines $D$ contained in $Z$ so, by lemma 3.2, it is enough to show
that $S_0$ is adjusted in $\sbm P^n$ in degree d. We will apply lemma 3.1 to
$S_0$.

The trace $S_0\upp$ of $S_0$ on $H$, is the generic union of $u+v-s$
first infinitesimal neighbourhoods of closed points in $H$ with $$a\up=
a-t+b\geq(n-1)(n-2)/2$$ closed points and an $(s-2)$-complex. The scheme
$S_0\upp$ is numerically adjusted in degree $d=3$ in $H= \sbm P^{n-1}$ by
$(\alpha)$ of lemma 6.3.4. By specialising closed points to first
infinitesimal neighbourhoods we can suppose that 
\begin{equation} ((n-1)(n-2)/2) + (n-1)  \geq   a\up  \geq  (n-1)(n-2)/2
\end{equation}
Hence by (15) we conclude that $S_0\upp$ admits a specialisation 
to a subscheme of $H$ satisfying the conditions of the proposition in
dimension $n-1$. Using the induction hypothesis we conclude that
$S_0\upp$ is adjusted in $H$ in degree $d=3$.

By lemma 3.1, the residual $S_0\up$ of $So$ with respect to H is
numerically adjusted in $\sbm P^n$ in degree $d=2$. It is the generic union
of the first infinitesimal neighbourhoods of $1\leq s \leq (n+1)/2$ closed 
points of $\sbm P^n$ and the $s(s-1)/2$ lines which join them two by two,
with an $(l-2)$-complex, $u+v-s+b$ closed points in $H$ and $t$ finite,
degree two, closed subschemes of $\sbm P^n$ with support in $H$, but
transverse to $H$. Let $T$ be the specialisation of $S_0\up$ obtained by
specialising all but one of the infinitesimal neighbourhoods to have support
in $H$ and specialising all the degree two finite closed subschemes into
$H$. By [AH2] (4.9) the trace $T\upp$ on $H$ is adjusted in $H$ in degree
$d=2$ and it is clear that the residual $T\up$ of $T$ with respect to $H$
(which is the union of a first infinitesimal neighbourhood
of a closed point of $\sbm P^n$ with $(s-1)$ lines through the point) is
adjusted in degree one in $\sbm P^n$. Applying lemma 3.1 to $T$ we conclude
that $S_0\up$ is adjusted in degree $d=2$ in $\sbm P^n$. 

Now applying lemma 3.1 to $S_0$ we conclude that $S_0$ is 
adjusted in degree $d=3$ in $P^n$.$\quad\quad\Box$
	
\noindent{\bf References}

\putref{[A]}{J. Alexander:
{\em Singularit\'es imposables en position g\'en\'erale aux  hypersurfaces de
$\sbm P^n$}. Compos. Math. 68 (1988) 305-354.}

\putref{[AH1]}{J. Alexander- A. Hirschowitz: {\em Un lemme de Horace
diff\'erentiel: application aux singularit\'es des hyperquartiques de $\sbm
P^5$}. Journal of Alg. Geom. vol 1. N$^o$ 3 (1992) 411-426}

\putref{[AH2]}{ J. Alexander-
A. Hirschowitz: {\em Polynomial Interpolation in Several 	Variables.}
Journal of Alg. Geom.  vol 4. N$^o$ 4 (1995) 201-222 }

\putref{[AH3]} {J. Alexander-
A. Hirschowitz: {\em La m\'ethode d'Horace \'eclat\'ee: 	application aux
singularit\'es des hyperquartiques.} Invent. 	Math. 107 (1992) 585-602} 

\putref{[CH]}{C. Ciliberto and A. Hirschowitz: {\em Hypercubiques de $\sbm
P^4$ avec sept points singuliers g\'en\'eriquaes}. CRAS 313 (1991) 135-137}

\putref{[H1] }{A. Hirschowitz: La
m\'ethode d'Horace pour l'interpolation \`a plusieurs 	variables. Manuscr.
Math.  50 (1985) 337-388.} 

\putref{[H2]}{A. Hirschowitz: Rank techniques and jump
stratifications. In: Vector  bundles on algebraic varieties, Proceedings
Bombay  1984, 	Oxford  Univ. Press (1987) 159-205.} 

\putref{[H3] }{A. Hirschowitz:
Une conjecture pour la cohomologie des diviseurs sur 	les surfaces
rationnelles g\'en\'eriques. J. fŸr die r. und ang. Math.
	397 (1989) 208-213.}

\begin{minipage}[t]{.4\linewidth}{James
Alexander\\
Universit\'e d'Angers\\
jea@tonton.univ-angers.fr}\end{minipage}\hfill
\begin{minipage}[t]{.4\linewidth}{Andr\'e
Hirschowitz\\
Universit\'e de Nice\\
ah@math.unice.fr}\end{minipage}

\end{document}